\documentclass[12pt,draft]{amsart}
\usepackage{amsmath,amsthm,latexsym,amscd,amsbsy,amssymb,amsfonts,
euscript}

\renewcommand{\AE}{\mathop{\rm AE}\nolimits}
\newcommand{\ANE}{\mathop{\rm ANE}\nolimits}
\newcommand{\AnE}{\mathop{\rm A(N)E}\nolimits}

\newcommand{\dimm}{\mathop{\rm dim}\nolimits}

\newcommand{\cdim}{\mathop{\rm c-dim}\nolimits}

\newcommand{\U}{{\EuScript U}}

\newcommand{\ed}{\mathop{\rm ed}\nolimits}

\newcommand{\cit}[1]{{\rm {(\cite{#1})}.}}
\newcommand{\cov}{\mathop{\rm cov}\nolimits}

\newcommand{\al}{\alpha}

\newcommand{\z}{\zeta}

\newcommand{\tor}{\mathop{\rm Tor}\nolimits}

\newcommand{\Z}{\mathbb Z}
\newcommand{\ZZ}[1]{\mathbb Z _{#1}}
\newcommand{\fnd}{\pi_n^{[L]}(S^n)}
\newcommand{\group}{\pi_n^{[L]}(X)}
\newcommand{\dd}[2]{\delta ^*_{#1,#2}}
\newcommand{\mesh}{\mathop{\rm mesh}\nolimits}
\newcommand{\idn}{\mathop{\rm id}\nolimits}
\newcommand{\nsph}{S^n_{[L]}}

\newtheorem{thm}{Theorem}[section]
\newtheorem{cor}[thm]{Corollary}
\newtheorem{lem}[thm]{Lemma}
\newtheorem{pro}[thm]{Proposition}

\theoremstyle{definition}
\newtheorem{dfn}{Definition}[section]
\theoremstyle{remark}
\newtheorem{rem}{Remark}[section]

\chardef\bslash=`\\ 

\makeatletter
\def\verbatim{\interlinepenalty\@M \@verbatim
  \leftskip\@totalleftmargin\advance\leftskip2pc
  \frenchspacing\@vobeyspaces \@xverbatim}
\makeatother
\numberwithin{equation}{section}

\begin{document}


\title[On L-homotopy groups]
{On $[L]$-homotopy groups}
\author{A.~V.~Karasev}
\address{Department of Mathematics and Statistics,
University of Saskatche\-wan,
McLean Hall, 106 Wiggins Road, Saskatoon, SK, S7N 5E6,
Canada}
\email{karasev@math.usask.ca}



\begin{abstract}
The paper is devoted to investigation of some properties of
$[L]$-homotopy groups. It is proved, in particular, that for any finite
$CW$-complex $L$, satisfying
double inequality $[S^n] < [L] \le [S^{n+1}]$,
$\fnd = \Z$. Here $[L]$ denotes extension type of complex $L$ and
$\group$ denotes $n$-th $[L]$-homotopy group of $X$.
\end{abstract}

\maketitle
\markboth{A.~V.~Karasev}{On $[L]$-homotopy groups}


\section{Introduction}

A new approach to dimension theory, based on notions of extension
types of complexes and extension dimension leads to
appearence of $[L]$-homotopy theory which, in turn, allows to introduce
$[L]$-homotopy groups (see \cite{ch}).
Perhaps the most natural problem related to $[L]$-homotopy groups is a
problem of computation. It is necessary to point out that $[L]$-homotopy
groups may differ from usual homotopy groups even for complexes.

More specifically the problem of computation can be stated as follows:
describe $[L]$-homotopy groups of a space $X$ in terms of usual homotopy
groups of $X$ and homotopy properties of complex $L$.

The first step on this way is apparently computation of $n$-th
$[L]$-homotopy group of $S^n$ for complex whose extension type lies
between extension types of $S^n$ and $S^{n+1}$.

In what follows we, in particular, perform this step.

\section{Preliminaries}

     Follow \cite{ch}, we introduce notions
of {\it extension types of complexes, extension dimension,
$[L]$-homotopy, $[L]$-homotopy groups} and other related notions.

We also state Dranishnikov's theorem, characterizing
extension properties of complex \cite{dr}.

All spaces are polish,
all complexes are countable finitely-dominated $CW$ complexes.

For spaces $X$ and $L$, the notation $L \in \AE (X)$ means, that every map
$f:A \to L$, defined on a closed subspace $A$ of $X$, admits an
extension $\bar f$ over $X$.

Let $L$ and $K$ be complexes. We say (see \cite{ch}) that $L \le K$ if
for each space $X$ from $L \in \AE (X)$ follows $K \in \AE (X)$. Equivalence
classes of complexes with respect to this relation are called {\it
extension types}. By $[L]$ we denote extension type of $L$.

\begin{dfn} \cit {ch} The extension dimension of a space $X$ is
extension type $\ed (X)$ such that
$\ed (X) = \min \{ [L] : L \in \AE  (X) \}$.
\end{dfn}

Observe, that if $[L] \le [S^n]$ and $\ed (X) \le [L]$, then $\dimm X
\le n$.

Now we can give the following

\begin{dfn}\cite{ch}
We say that a space $X$ is {\it an absolute (neighbourhood) extensor modulo}
$L$ (shortly $X$ is  $\AnE ([L])$)
and write $X\in \AnE ([L])$ if $X\in \AnE (Y)$ for each space
$Y$  with $\ed (X) \le [L]$.
\end{dfn}

Definition of $[L]$-homotopy and $[L]$-homotopy equivalence \cite{ch}
are essential for our consideration:

\begin{dfn}
Two maps $f_0$, $f_1: X \to Y$ are said to be $[L]$-homotopic (notation:
$f_0 \stackrel{[L]}{\simeq} f_1$) if for any map
$h: Z \to X \times [0,1]$, where $Z$ is a space with $\ed (Z) \le [L]$,
the composition $(f_0 \oplus f_1) h|_{h^{-1}(X \times \{ 0, 1 \} )} :
h^{-1} (X \times \{ 0,1 \}) \to Y$ admits an extension $H: Z \to Y$.
\end{dfn}

\begin{dfn}
A map $f:X \to Y$ is said to be $[L]$-homotopy equivalence if there is a
map $g: Y \to X$ such that the compositions $gf$ and $fg$ are
$[L]$-homotopic to $\idn _X$ and $\idn _Y$ respectively.
\end{dfn}

Let us observe (see \cite{ch}) that $\ANE ([L])$-spaces have the following
$[L]$-homotopy extension property.

\begin{pro} \label{homext}
Let $[L]$ be a finitely dominated complex and $X$ be a Polish $\ANE
([L])$-space. Suppose that $A$ is closed in a space $B$ with
$\ed (B) \le [L]$. If maps $f, g : A \to X$ are $[L]$-homotopic and
$f$ admits an extension $F : B \to X$ then $g$ also admits an
extension $G : B \to X$, and it may be assumed that $F$ is
$[L]$-homotopic to $G$.
\end{pro}

To provide an important example of $[L]$-homotopy equivalence we need
to introduce the class of approximately $[L]$-soft maps.

\begin{dfn}\cite{ch}
A map $f: X \to Y$ is said to be approximately $[L]$-soft, if for each
space $Z$ with $\ed (Z) \le [L]$, for each closed subset $A\subset Z$,
for an open cover $\U \in \cov (Y)$, and for any two maps $g: A \to X$
and $h: Z \to Y$ such that $fg = h| _A$ there is a map $k: Z \to X$
satisfying  condition $k| _A = g$ and the composition $fk$ is $\U$-close
to $h$.
\end{dfn}

\begin{pro}\cite{ch} \label{hmteq}
Let $f: X \to Y$ be a map between $\ANE ([L])$-compacta and
$\ed (Y) \le [L]$. If $f$ is approximately $[L]$-soft then
$f$ is a $[L]$-homotopy equivalence.
\end{pro}

In order to define $[L]$-homotopy groups it is necessary to consider
an {\it $n$-th $[L]$-sphere $\nsph$} \cite{ch},
namely, an $[L]$-dimensional $\ANE
([L])$ - compactum admitting an approximately $[L]$-soft map onto $S^n$.
It can be shown that all possible choices of an $[L]$-sphere $\nsph$ are
$[L]$-homotopy equivalent. This remark, coupled with the following
proposition, allows us to consider for every finite complex $L$,
every $n \ge 1$ and for
any space $X$, the set $\group = [\nsph, X] _{[L]}$ endowed with natural
group structure (see \cite{ch} for details).

\begin{thm} \cite{ch}
Let $L$ be a finitely dominated complex and $X$
be a finite polyhedron or a compact Hilbert cube manifold.
Then there exist a $[L]$-universal $\ANE ([L])$ compactum $\mu ^{[L]}
_X$ with $\ed (\mu ^{[L]} _X ) = [L]$ and
 an $[L]$-invertible and approximately $[L]$-soft map
$f^{[L]}_X : \mu ^{[L]}_X \to X$.
\end{thm}

The following theorem is essential for our consideration.

\begin{thm} \label{D1}
Let $L$ be simply-connected $CW$-complex, $X$ be finite-dimensional
compactum. Then $L \in \AE (X)$ iff $\cdim _{H_i(L)} X \le i$ for any $i$.
\end{thm}

From the proof of Theorem \ref{D1} one can conclude that the following
theorem also holds:

\begin{thm} \label{D2}
Let $L$ be a $CW$-complex (not necessary\\ simply-connected). Then for any
finite-dimensional compactum $X$ from $L \in \AE (X)$ follows that
$\cdim _{H_i(L)} X \le i$ for any $i$.
\end{thm}

\section{Cohomological properties of $L$}

In this section we will investigate some cohomological properties of
complexes $L$ satisfying condition $[L] \le S^n$ for some $n$.
To establish these properties let us first formulate the following

\begin{pro} \label{Sp} \cite{sp}
Let $(X,A)$ be a topological pair, such that $H_q (X,A)$ is finitely
generated for any $q$. Then free submodules of $H^q (X,A)$ and
$H_q (X,A)$ are isomorphic and torsion submodules of $H^q (X,A)$ and
$H_{q-1} (X,A)$ are isomorphic.
\end{pro}

Now we use Theorem \ref{D2} to obtain the following lemma.

\begin{lem} \label{torh}
Let $L$ be finite $CW$ complex such that $[L] \le [S^{n+1}]$  and $n$ is
minimal with this property. Then for any $q \le n$ $H_q (L)$ is torsion
group.
\end{lem}
\begin{proof}
Suppose that there exists $q \le n$ such that $H^q (L) = \Z \oplus G$.
To get a contradiction let us show that $[L] \le [S^q]$. Consider $X$
such that $L \in \AE (X)$. Observe, that $X$ is finite-dimensional since
$[L] \le [S^{n+1}]$ by our assumption.

Denote $H = H_q (L)$. By Theorem \ref{D2} we have $\cdim _H X \le q$.
Hence, for any closed subset $A \subseteq X$  we have
$H^{q+1} (X,A;H) = \{ 0 \}$. From the other hand, univeral coefficients
formula implies that\\
$H^{q+1} (X,A) \approx H^{q+1} (X,A) \otimes H \oplus \tor (H^{q+2} (X,
A),H)$.

Hence, $H^{q+1} (X,A) \otimes H = \{ 0 \}$. Observe, however, that by our
assumtion we have $H^{q+1} (X,A) \otimes H = H^{q+1} \otimes ( \Z
\oplus G) = H^{q+1} (X,A) \oplus (H^{q+1} (X,A) \otimes G)$.
Therefore, $H^{q+1} (X,A) = 0$.

From the last fact we conclude that $\cdim X \le q$ and therefore since
$X$ is finite-dimensional, $\dimm X \le q$ which iplies $S^q \in \AE
(X)$.
\end{proof}

From this lemma and Proposition \ref{Sp} we obtain

\begin{cor} \label{torch} In the same assumptions $H^{q} (L)$ is
torsion group for any $q \le n$.
\end{cor}

The following fact is essential for constraction of compacts with some
specific properties which we are going to construct further.

\begin{lem}\label{acycl}
Let $L$ be as in previous lemma. For any $m$ there exists $p \ge m$ such
that $H^q(L; \ZZ{p} ) = \{ 0 \}$ for any $q \le n$.
\end{lem}
\begin{proof}
From Corollary \ref{torch} we can conclude that $H^{q} (L) =
\bigoplus\limits_{i=1}^{l_k} \ZZ{m_{qi}}$ for any $q \le n$. Additionally,
let $\tor H^{n+1}(L) = \bigoplus\limits_{i=1}^{l_{n+1}} \ZZ{m_{(n+1)i}}$

For any $m$ consider $p \ge m$ such that $(p,m_{ki})$ for every
$k = 1 \ldots n+1$ and $i = 1 \ldots l_k$.
Universal coefficients formula implies that $H^{q} (L; \ZZ{p} ) = \{ 0
\}$ for every $k \le n$.
\end{proof}

Finally let us proof the following

\begin{lem}\label{ext}
Let $X$ be a metrizable compactum, $A$ be a closed subset of $X$. Consider
a map $f: A \to S^n$. If there exists extension $\bar f : X \to S^n$ then
for any $k$ we have $\dd{X}{A} (f^* ( \z )) = 0$ in group $H^{n+1} (X,A;
\ZZ{k})$, where $\z$ is generator in $H^n (S^n, \ZZ{k})$.
\end{lem}
\begin{proof}
Let $\bar f$ be an extension of $f$.
Commutativity of the following diagram implies assertion of lemma:

$$
\begin{CD}
H^n (A; \ZZ {k}) @>\dd{X}{A}>> H^{n+1} (X,A; \ZZ {k})\\
@AA{\bar{f^*} = f^*}A   @AA{\bar{f^*}}A\\
H^n (S^n; \ZZ {k}) @>\dd{S^n}{S^n}>> H^{n+1} (S^n, S^n; \ZZ {k}) = \{ 0 \}
\end{CD}
$$

\end{proof}

\section{Some properties of [L]-homotopy groups}

In this section we will investigate some properties of $[L]$-homotopy
groups.

From this point and up to the end of the text we consider finite complex
$L$ such that $[S^n] < [L] \le [S^{n+1}]$ for some fixed $n$.

\begin{rem} \label{grup}
Let us observe that for such complexes $\nsph$ is $[L]$-homotopic
equivalent to $S^{n}$ (see Proposition \ref{hmteq}).
Therefore for any $X$ $\group$ is isomorphic to $G = \pi _n (S^n)/N([L])$
where $N([L])$ denotes the relation of $[L]$-homotopic equivalence
between elements of $\pi_n (S^n)$.
\end{rem}

From this observation one can easely obtain the following fact.

\begin{pro}\label{grp}
For $\pi_n^{[L]}(S^n)$ there are three variants:
$\pi_n^{[L]} (S^n) = \Z$,  $\pi_n^{[L]} (S^n) = \ZZ{m}$ for
some integer $m$ or this group is trivial.
\end{pro}

Let us characterize the hypothetical equality $\fnd = \ZZ{m}$
in terms of extensions of maps.

\begin{pro}\label{hext}
If $\fnd = \ZZ{m}$ then for any $X$ such that  $\ed (X) \le [L]$,
for any closed subset $A$ of $X$ and for any map $f:A \to S^n$,
there exists
extension $\bar h : X \to S^m$ of composition $h = z_m f$, where $z_m :
S^n \to S^n$ is a map having degree $m$.
\end{pro}
\begin{proof}
Suppose, that $\fnd = \ZZ{m}$. Then from Remark \ref{grup}
and since $ [z_m] = m[ \idn _{S^n}] = [*]$ (where $[f]$ denotes homotopic
class of $f$) we conclude that $z_m : S^n \to S^n$ is $[L]$-homotopic to
constant map. Let us show that $h = z_m f : A \to S^n$ is also
$[L]$-homotopic to constant map. This fact will prove our statement.
Indeed, by our assumption $\ed (X) \le [L]$ and $S^n \in ANE$ and
therefore we can apply Proposition \ref{homext}.

Consider $Z$ such that $\ed (Z) \le [L]$ and a map $H : Z \to A \times
I$, where $I = [0,1]$. Pick a point $s \in S^n$.
Let $f_0 = z_m f$, $f_1 \equiv s$ -- constant map considered as $f_i : A
\times \{ i \} \to S^n$, $i = 0,1$.

Define $F: A \times I \to S^n \times I$ as follows: $F (a,t) = (f(a),t)$
for each $a \in A$ and $t \in I$.
Let $f_0' \equiv z_m$ and $f_1' \equiv s$ considered as $f_i' : S^n
\times \{ i \} \to S^n$, $i = 0,1$.

Consider a composition $G = F H : Z \to S^n \times I$.
By our assumption $f_0'$ is $[L]$-homotopic to $f_1'$. Therefore
a map $g: G^{-1} (S^n \times \{ 0 \} \bigcup S^n \times \{ 1 \}) \to
S^n$, defined as $g| _{G^{-1} (S^n \times \{ i \}) } = f_i' G$ for $i =
0,1$, can be extended over $Z$. From the other hand we have
$G^{-1} (S^n \times \{ i \}) \equiv H^{-1} (A \times \{ i \})$ and
$g| _{G^{-1} (S^n \times \{ i \}) } = f_i' f H = f_i$ for $i =
0,1$. This remark completes the proof.
\end{proof}

Now consider a special case of complex having a form
$S^n < L = K_s \vee K \le S^{n+1}$,
where $K_s$ is a complex obtained by attaching to $S^n$ a $(n+
1)$-dimensional cell using a map of degree $s$.

\begin{pro}
Let $[\al ] \in \pi _n (X)$ be an element of order $s$. Then $\al$ is
$[L]$-homotopy to constant map.
\end{pro}
\begin{proof}
Observe that simillar to proof of Proposition \ref{hext} it is enough to
show that for every $Z$ with $\ed (Z) \le [L]$, for every closed
subspace $A$ of $Z$ and for any map $f : Z \to S^n$ a composition $\al f
: A \to X$ can be extended over $Z$.

Let $g : S^n \to K_s ^{(n)}$ be an embedding (by $M ^{(n)}$ we denote
$n$-dimensional skeleton of complex $M$) and $r : L \to K_s$ be a
retraction.

Since $\ed (Z) \le [L]$, a composition $gf$ has an extension $F : Z \to
L$. Let $F' = rF$ and $\al '$ be a map $\al$ considered as a map $\al '
: K_s ^{(n)} \to X$.
Observe that $\al' F'$ is a necessary extension of $\al f$.
\end{proof}

\section{Computation of $\fnd$}

In this section we will prove that $\fnd = \Z$.

Suppose the oppsite, i.e. $\fnd = \ZZ{m}$ (we use Proposition \ref{grp};
the same arguments can be used to prove that $\fnd$ is non-trivial).

To get a contradiction we need to obtain a compact with special
extension properties. We will use a construction of \cite{drp}

Let us recall the following definition.

\begin{dfn} \cite{drp}
Inverse sequence $S = \{ X_i, p ^{i+1}_i : i \in \omega \}$ consisting
of metrizable compacta is said to be $L$-resolvable if for any $i$, $A
\subseteq X_i$ - closed subspace of $X_i$ and any map $f:A \to L$ there
exists $k \le i$ such that composition $fp ^k _i : (p ^k _i)^{-1}A \to
L$ can be extended over $X_k$.
\end{dfn}

The following lemma (see \cite{drp}) expresses an important property of
$[L]$-resolvable inverse sequences.

\begin{lem}
Suppose that $L$ is a countable complex and that $X$ is a compactum such
that $X = lim S$ where $S = { (X_i, \lambda _i), q_i^{i+1}}$ is a
$L$-resolvable inverse system of compact polyhedra $X_i$ with
triangulations $\lambda _i$ such that $\mesh \{ \lambda _i \} \to 0$.
Then $L \in \AE (X)$
\end{lem}

Let us recall that in \cite{drp} inverse sequence
$S = \{ (X_i, \tau _i), p^{i+1}_i \}$ was constructed
such that $X_i$ is compact polyhedron with fixed
triangulation $\tau _i$, $X_0 = S^{n+1}$,
$\mesh \tau _i \to 0$, $S$ is $[L]$-resolvable and
for any $x \in X_i$ we have $(p^{i+1}_i)^{-1}x \simeq L$ or $*$.

It is easy to see
that using the same  construction one can
obtain inverse sequence $S = \{ (X_i, \tau _i), p^{i+1}_i \}$
having the same properties with exeption of
$X_0 = D^{n+1}$ where $D^{n+1}$ is $n+1$-dimensional disk.

Let $X = \lim S$. Observe, that $\ed (X) \le [L]$. Let $p_0 : X \to
D^{n+1}$ be a limit projection.

Pick $p \ge m+1$ which Lemma \ref{acycl} provides us with.
By Vietoris-Begle theorem (see \cite{sp}) and our choice of p,
for every $i$ and every $X'_i \subseteq X_i$ a
homomorphism $(p^{i+1}_i)^* : H^{k} (X'_i; \ZZ {p} )
\to H^{k} ((p^{i+1}_i)^{-1}X'_i ; \ZZ{p} )$ is isomorphism for $k \le n$
and monomorphism for $k = n+1$.

Therefore for each $D' \subseteq X_0 = D^{n+1}$ homomorphism
$p_0^* : H^{k} (D'; \ZZ{p} ) \to H^{k} ((p_0)^{-1} D' ; \ZZ{p}) $ is
isomorphism for $k \le n$ and monomorphism for $k = n+1$. In particular,
$H^{n} (X ; \ZZ{p}) = \{ 0 \}$ since $X_0 = D^{n+1}$ has trivial
cohomology groups.

Let $A = (p_0)^{-1} S^n$
and $\z \in H^n(S^n; \ZZ{p}) \approx \ZZ{p}$  be a generator.

Since $p_0^* : H^n (S^n; \ZZ{p} ) \to H^n (A; \ZZ{p} )$ is isomorphism,
$p_0^* (\z )$ is generator in $H^n (A, \ZZ{p}) \approx \ZZ{p} $.
In particular, $p_0^* (\z )$ is element of
order $p$.

From exact sequence of pair $(X,A)$

$$
\begin{CD}
\ldots\to H^n (X; \ZZ{p}) = \{ 0 \} @>i_{X,A}>> H^n (A; \ZZ {p})
@>\dd{X}{A}>> H^{n+1} (X,A; \ZZ{p}) \to \ldots \\
\end{CD}
$$

we conclude that $\dd{X}{A}$ is monomorphism and hence
$\dd{X}{A}(p_0^*(\z )) \in H^{n+1} (X,A; \ZZ{p})$ is element
of order $p$.

Consider now a composition $h = z_m p_0$. By our assumption
this map can be extended over $X$ (see Proposition \ref{hext}).
This fact coupled with  Lemma \ref{ext} implies that
$\dd{X}{A}(h^*(\z )) = {0}$ in $H^{n+1} (X,A; \ZZ{p})$.
But $\dd{X}{A}(h^*(\z )) = m\dd{X}{A}(p_0^*(\z ))$.
We arrive to a contradiction which shows that

\begin{thm}
Let $L$ be a complex such that $[S^n] < [L] \le [S^{n+1}]$. Then
$\fnd \approx \Z$.
\end{thm}

The author is greatfull to A.~C.~Chigogidze for usefull discussions.



\begin{thebibliography}{99}

\bibitem{ch}
A.~Chigogidze, {\em Infinite dimensional topology and shape theory},
to appear in: "Handbook of Geometric Topology" edited by R.~Daverman and
R.~B.~Sher, North Holland, Amsterdam, 1999.

\bibitem{dr}
A.~N.~Dranishnikov, {\em Extension of mappings into $CW$-complexes},
Math. USSR Sbornik {\bf 74} (1993), 47-56.

\bibitem{drp}
A.~N.~Dranishnikov and D.~Repov\v{s},
{\em Cohomological dimension with respect to perfect groups},
Topology Appl. {\bf 74} (1996), 123-140.

\bibitem{sp}
E.~H.~Spanier,
{\em Algebraic topology}, McGraw-Hill, New York, 1966.

\end{thebibliography}
\end{document}